\documentclass[a4paper,12pt,leqno,twoside]{article}

\usepackage{amsmath, amsthm, amssymb}

\numberwithin{equation}{section}

\def\bN{\mathbb{N}}
\def\bZ{\mathbb{Z}}

\def\bP{\mathbb{P}}

\def\P{\mathcal{P}}
\def\d{\mathrm{d}}
\def\PP{\mathrm{P}}
\def\Prob{\mathrm{Prob}}

\def\({ \left( }
\def\){ \right)}

\def\Gessel{%
 \begin{subarray}{c}
 \lambda \in \P \\
 \lambda_1 \leq h 
 \end{subarray}
}

\theoremstyle{plain}

\newtheorem{thm}{Theorem} 

\newtheorem{lem}{Lemma} 

\theoremstyle{definition}

\newtheorem{remark}{Remark} 

\title{A Scaling Limit for $t$-Schur Measures}
\author{\textsc{Sho Matsumoto}}

\date{\today}

\begin{document}
\setlength{\baselineskip}{16pt}
\maketitle

\begin{abstract}
We introduce a new measure on partitions.
We assign
to each partition $\lambda$ a probability $S_{\lambda}(x;t) s_{\lambda}(y) / Z_{t}$
where $s_{\lambda}$ is the Schur function,
$S_{\lambda} (x ;t)$ is a generalization of the Schur function defined in \cite{Macdonald} and 
$Z_{t}$ is a normalization constant.   
This measure, which we call the $t$-Schur measure, is a generalization of the Schur measure \cite{Okounkov}
and is closely related to
the shifted Schur measure studied by Tracy and Widom \cite{TracyWidom200?} for a combinatorial viewpoint.

We prove that
a limit distribution of the length of the first row of a partition with respect to
$t$-Schur measures is given by
the Tracy-Widom distribution, i.e.,
the limit distribution of the largest eigenvalue suitably centered and normalized 
in GUE.
\end{abstract}

 
%
%
\section{Introduction}
%

Let $\P$ be the set of all partitions $\lambda$ and $s_{\lambda}$ the Schur function (see \cite{Macdonald})
with variables $x=(x_1, x_2, \dots )$ or $y=(y_1, y_2, \dots )$.
The Schur measure introduced in \cite{Okounkov} is a probability measure on $\P$ defined by 
 \begin{equation} \label{Schurmeasure}
   \PP_{\mathrm{Schur}}(\{ \lambda \}) := \frac{1}{Z_0} s_{\lambda} (x) s_{\lambda} (y),
 \end{equation}
where the normalization constant $Z_0$ is determined by the Cauchy identity
 \begin{equation}\label{CauchySchur}
   Z_0 := \sum_{\lambda \in \P} s_{\lambda}(x) s_{\lambda}(y) = \prod_{i,j=1}^{\infty} \frac{1}{1 - x_i y_j}. 
 \end{equation}

We consider a certain specialization of this measure.
Let $\alpha$ be a real number such that $0 < \alpha <1$.
We put $x_i = \alpha$ and $y_j = \alpha$ for $1 \leq i \leq m$ and $1 \leq j \leq n$,
and let the rest be zero.
Fix $\tau=m/n$. 
This is called the \textit{$\alpha$-specialization}.

Johansson \cite{Johansson1} (see also \cite{Johansson2}, \cite{Johansson3})
showed that
when $n \to \infty$ 
a distribution of the length of the first row $\lambda_1$ of a partition $\lambda$
with respect to the $\alpha$-specialized Schur measure 
converges to
the Tracy-Widom distribution \cite{TracyWidom1},
which is the limit distribution of the largest eigenvalue suitably centered and normalized
in the Gaussian Unitary Ensemble (GUE).

The Tracy-Widom distribution $F_2(s)$ is explicitly expressed as
      $$
       F_2(s) := \sum_{k=0}^{\infty} \frac{(-1)^k}{k !} \int_{[s, \infty )^k} \det (K_{\mathrm{Airy}}(x_i, x_j))_{i,j=1}^k
       \d x_1 \dots \d x_k,
      $$
    where $K_{\mathrm{Airy}}(x,y)$ denotes the Airy kernel given by
      $$ 
       K_{ \mathrm{Airy} } (x,y) = \int_0^{\infty} \mathrm{Ai}(x+z) \mathrm{Ai} (z+y) \d z
      $$  
    and $\mathrm{Ai}(x)$ denotes the Airy function given by 
      $$
       \mathrm{Ai}(x)= \frac{1}{2 \pi \sqrt{-1}} \int_{\infty e^{-\pi \sqrt{-1}/3}}^{\infty e^{\pi \sqrt{-1}/3}}        
        e^{z^3/3 -x z} \d z.
      $$
On the other hand, 
Tracy and Widom \cite{TracyWidom200?} studied an analogue of the Schur measure,
which they call the shifted Schur measure,
and proved that a limit distribution of $\lambda_1$ with respect to the $\alpha$-specialized shifted Schur measure
is also given by the Tracy-Widom distribution.

In this paper,
we introduce a generalization of the Schur measure.
In order to define such a new measure,
let us recall the symmetric functions $e_n (x;t)$ with parameter $t$,
given by the generating function
 \begin{equation} \label{E}
    E_{x,t} (z) := \prod_{i=1}^{\infty} \frac{1+ x_i z}{1 + t x_i z} = \sum_{n=0}^{\infty} e_n (x;t) z^n.
 \end{equation}
Then
the (generalized) Schur function is given by
\begin{equation} \label{S}
    S_{\lambda} (x;t) := \det( e_{\lambda'_i -i +j}(x;t)),
 \end{equation}
where the partition $\lambda'$ is the conjugate of a partition $\lambda$, i.e.,
$\lambda'_i$ is the length of the $i$-th column of $\lambda$. 
These functions satisfy the so-called Cauchy identity (see \cite{Macdonald})
\begin{equation} \label{Cauchy}
   Z_t := \sum_{\lambda \in \P} S_{\lambda}(x;t) s_{\lambda} (y) = \prod_{i,j=1}^{\infty} \frac{1 -t x_i y_j}{1-x_i y_j}.
\end{equation}
In particular,
when $t=0$, $E_x(z)= E_{x,0}(z)$ is 
the generating function of elementary symmetric functions 
$e_n(x)=e_n(x;0) = \sum_{i_1<i_2< \dots <i_n} x_{i_1} x_{i_2} \dots x_{i_n}$.
Since $s_{\lambda}=\det(e_{\lambda'_i-i+j})$ (the dual version of the Jacobi-Trudi identity),
we notice that the identity \eqref{Cauchy} becomes \eqref{CauchySchur} when $t=0$.

By means of the identity \eqref{Cauchy},
we may define a probability measure on partitions $\lambda$ by
\begin{equation} \label{t-Schurmeasure}
    \PP_t (\{ \lambda \}) := \frac{1}{Z_t} S_{\lambda} (x;t) s_{\lambda} (y).
\end{equation}
We call this measure a \textit{$t$-Schur measure}. 
It reduces to the Schur measure at $t=0$.

Our main result is as follows.
Denote by $\PP_{\sigma}$ the $\alpha$-specialized $t$-Schur measure,
where $\sigma =  (m,n, \alpha, t)$ is the associated set of parameters of the measure.
\noindent
 \begin{thm} \label{mainthm}
    Suppose $-\infty < t \leq 0$.
    Then there exist positive constants 
    $c_1 =c_1 (\alpha, \tau, t)$ and $c_2 = c_2 (\alpha, \tau, t)$ such that
     $$
      \lim_{n \to \infty} \PP_{\sigma} \( \frac{ \lambda_1 -c_1 n}{c_2 n^{1/3}} <  s \) = F_2 (s).
     $$
 \end{thm}   
This theorem shows that the fluctuations in $\lambda_1$ are independent of the parameter $t$.
The assumption that $t$ is non-positive is required 
since the right hand side of \eqref{t-Schurmeasure} must be non-negative after making the $\alpha$-specialization.

In the case where $t=0$, this theorem gives the result due to Johansson \cite{Johansson1}. 
In Remark 2, we shall give the explicit expressions of $c_1 (\alpha, \tau, 0)$ and $c_2 (\alpha, \tau, 0)$ 
and explain the connection to the result in \cite{Johansson1}.
Note also that Theorem 1 does not imply the result of \cite{TracyWidom200?}
(see Remark 1). 
The proof of Theorem 1 will be given using the method of Tracy-Widom \cite{TracyWidom200?}.
The key of the proof is the determinantal expression of
$S_{\lambda}(x;t)$.

Further
we give a combinatorial interpretation of the $t$-Schur measure.
Namely, if we denote by $\bP$ an ordered set $\{ 1' < 1 < 2' < 2 < 3' <3 < \dots \}$,
then by virtue of the Robinson-Schensted-Knuth (RSK) correspondence 
between matrices with entries in $\bP \cup \{ 0 \}$ and pairs of tableaux of the same shape $\lambda=(\lambda_1, \lambda_2, \dots)$,
we see that the $t$-Schur measure corresponds to 
a measure (depending on $t$) on $\bP$-matrices.
According to this correspondence,
$\lambda_1$ corresponds to 
the length of the longest increasing subsequence in the biword $w_A$ associated with a $\bP$-matrix $A$.
Using the RSK correspondence and the shifted RSK correspondence,
we find that
this measure on $\bP$-matrices at $t=0$ and $t=-1$, respectively, corresponds to
the (original) Schur measure and
the shifted Schur measure,
respectively
(see \cite{Johansson1} and \cite{TracyWidom200?}).

%
%
\section{Schur functions, marked tableaux and the RSK correspondence}
%
%

In this section,
we summarize basic properties of Schur functions and marked tableaux
for providing a combinatorial interpretation of the $t$-Schur measure
(see \cite{Macdonald} and \cite{Sagan2} for details).

We denote the Young diagram of a partition $\lambda$ by the same symbol $\lambda$.
Let $\bN$ be the set of all positive integers 
and $\bP$ the totally ordered alphabet $\{ 1' < 1 < 2' < 2 < \dots \}$.
The symbols $1',2',3', \dots$ or $1, 2, 3, \dots$ are said to be \textit{marked} or \textit{unmarked}, respectively.
When it is not necessary to distinguish a marked element $k'$ from the unmarked one $k$, 
we write it by $|k|$.
A \textit{marked tableau} $T$ of shape $\lambda$ is 
an assignment of elements of $\bP$ to the squares of the Young diagram $\lambda$
satisfying the two conditions.
\begin{description}
\item[T1] The entries in $T$ are weakly increasing along each row and down each column.
\item[T2] For each $k \geq 1$, each row contains at most one marked $k'$ and 
         each column contains at most one unmarked $k$.
\end{description}   
The condition \textbf{T2} says that for each $k \geq 1$
the set of squares labelled by $k$ (resp. $k'$) is a horizontal (resp. vertical) strip.

For example,
 $$
  \begin{matrix}
   1' & 1 & 1 & 2' & 3' & 3 \\
   1' & 2 &   &    &    & \\
   3  & 3 &   &    &    &
  \end{matrix}
 $$
is a marked tableau of shape (6, 2, 2).

To each marked tableau $T$,
we associate a monomial
$ x^{T} = \prod_{i \geq 1} x_i^{m_i(T)}$,
where $m_i(T)$ is the number of times that $|i|$ appears in $T$. 
In the example above,
we have
$x^{T}= x_1^4 x_2^2 x_3^4$.

By the definition of $S_{\lambda} (x;t)$ and Chapter I, {\S}5, Example 23 in \cite{Macdonald},
it follows that
 \begin{equation} \label{generate}
  S_{\lambda} (x;t) = \sum_{T} (-t)^{\mathrm{mark}(T)} x^T,
 \end{equation}
where the sum runs over all marked tableaux $T$ of shape $\lambda$.
Here $\mathrm{mark} (T)$ is the number of marked entries in $T$.
In particular, we have
 $$
  s_{\lambda} (x) = S_{\lambda} (x;0) = \sum_{T} x^T,
 $$
where the sum runs over all marked tableaux which have no marked entries 
(i.e., all semi-standard tableaux) of shape $\lambda$.

We next explain the RSK correspondence between $\bP$-matrices and pairs of tableaux
(see \cite{Knuth}, \cite{Sagan1} and \cite{HoffmanHumphreys}).
Here $\bP$-matrix stands for the matrix whose entries are in $\bP_0 = \bP\cup \{0\}$.
To each $\bP$-matrix $A=(a_{ij})$ we associate a \textit{biword} $w_A$ as follows.
For $i,j \geq 1$,
the pair $\begin{pmatrix} i \\ j \end{pmatrix}$ is repeated $|a_{ij}|$ times in $w_A$, and
if $a_{ij}$ is marked, the lower entry $j$ of the first pair $\begin{pmatrix} i \\ j \end{pmatrix}$ appeared in $w_A$
is marked.
For example, 
 $$
   A=  \begin{pmatrix} 1' & 0 & 2 \\
                        2 & 1 & 2' \\
                       1' & 1' & 0 \end{pmatrix} \longmapsto
  w_A = \begin{pmatrix} 1 & 1 & 1 & 2 & 2 & 2 & 2 & 2 & 3 & 3 \\
                        1'& 3 & 3 & 1 & 1 & 2 & 3'& 3 & 1'& 2' \end{pmatrix}.
 $$ 
Observe that for a biword $w_A  = \begin{pmatrix} \beta_1 & \beta_2 & \dots & \beta_n \\
                                                  \alpha_1 & \alpha_2 & \dots & \alpha_n \end{pmatrix}$
the upper line $(\beta_1, \beta_2, \dots, \beta_n) \in \bN^n$ is a weakly increasing sequence.
Furthermore if $\beta_k = \beta_{k+1}$, then
$\alpha_k < \alpha_{k+1}$, or $\alpha_k$ and $\alpha_{k+1}$ are identical and unmarked.

Now we state a generalized RSK algorithm.
Let $S$ be a marked tableau and $\alpha$ an element in $\bP$. 
The procedure called an \textit{insertion} of $\alpha$ into $S$ is described as follows.
\begin{description}
  \item[I1] Set $R:=$ the first row of S.
  \item[I2] If $\alpha$ is unmarked, then
    \begin{description}
      \item[I2a] find the smallest element $\beta$ in $R$ greater than $\alpha$ and replace $\beta$ by $\alpha$ in $R$.
             (This operation is called the BUMP.)
      \item[I2b] set $\alpha := \beta$ and $R :=$ the next row down.  
    \end{description}
  \item[I3] If $\alpha$ is marked, then
    \begin{description}
      \item[I3a] find the smallest element $\beta$ in $R$ which is greater than or equal to $\alpha$ and replace $\beta$ by $\alpha$ in $R$.
                 (This is called the EQBUMP.)
      \item[I3b] set $\alpha := \beta$ and $R :=$ the next row down.
    \end{description}
  \item[I4] If $\alpha$ is unmarked and is greater than or equal to the rightmost element in $R$, or    
            if $\alpha$ is marked and greater than every element of $R$, then
            place $\alpha$ at the end of the row $R$ and stop.
\end{description} 
Write the result of inserting $\alpha$ into $S$ by $I_{\alpha} (S)$.

For a given biword $w_A =  \begin{pmatrix} \beta_1 & \beta_2 & \dots & \beta_n \\
                                                  \alpha_1 & \alpha_2 & \dots & \alpha_n \end{pmatrix}$,
we construct a sequence of pairs of a marked tableau and a semi-standard tableau as
 $$
  (S_0, U_0)=( \emptyset, \emptyset),\ (S_1, U_1), \ \dots, \ (S_n, U_n)=(S,U).
 $$
Assuming that a pair $(S_{k-1},U_{k-1})$ of the same shape is given.
Then we construct $(S_k,U_k)$ as follows.
A marked tableau $S_k$ is $I_{\alpha_k} (S_{k-1})$.
A semi-standard tableau $U_k$ is obtained by 
writing $\beta_k$ into the new cell of $U_k$  
created by inserting $\alpha_k$ to $S_{k-1}$.
We call $S=S_n$ a \textit{insertion tableau} and $U=U_n$ a \textit{recording tableau}.

For example,
for a biword $ w_A = \begin{pmatrix} 1 & 1 & 1 & 2 & 2 & 2 & 2 & 2 & 3 & 3 \\
                        1'& 3 & 3 & 1 & 1 & 2 & 3'& 3 & 1'& 2' \end{pmatrix}$, 
we obtain
$$ 
 (S,U)=\( \begin{matrix} 1'& 1 & 1 & 2'& 3'& 3 \\
                         1'& 2 &   &   &   &   \\
                         3 & 3 &   &   &   &   \end{matrix}, \quad
     \begin{matrix} 1 & 1 & 1 & 2 & 2 & 2 \\
                    2 & 2 &   &   &   &   \\
                    3 & 3 &   &   &   &   \end{matrix}\).
$$

The generalized RSK correspondence is then described as follows.

\noindent
\begin{thm}
 There is a bijection between $\bP$-matrices $A=(a_{ij})$ and  pairs $(S,U)$ of a marked tableau $S$ and
 an unmarked tableau $U$ of the same shape such that
 $\sum_{i} |a_{ij}|=s_j$ and $\sum_{j}|a_{ij}| = u_i$. 
 Here we put $s_k = m_k(S)$ and $u_k = m_k(U)$ for $k \geq 1$,
 and the number of marked entries in $A$ is equal to $\mathrm{mark}(S)$.
\end{thm}
\begin{proof}
The proof of this theorem can be done in a similar way to the original RSK correspondence 
between $\bN$-matrices and pairs of semi-standard tableaux, 
or of the shifted RSK correspondence between $\bP$-matrices and pairs of shifted marked tableaux
(See \cite{HoffmanHumphreys}, \cite{Sagan1}, \cite{Sagan2}). 
We omit the detail.
\end{proof}
We call the tableau $S$ (resp. $U$) to be of type $s=(s_1, s_2 , \dots)$ (resp. $u=(u_1, u_2, \dots)$).

An important property of this correspondence is 
its relationship to the length of the longest increasing subsequence in a biword $w_A$.
The {\it increasing subsequence} in $w_A$ is 
a weakly increasing subsequence in the lower line in $w_A$  such that a marked $k'$ appears at most one for each positive integer $k$. 
In the example above, 
$ (1'\  1 \  1 \ 2 \ 3' \ 3)$ is one of such increasing subsequences in $w_A$. 
Let $\ell (w_A)$ denote the length of the longest increasing subsequence in $w_A$.
Then we have the
 
\noindent
\begin{thm}
  If a $\bP$-matrix $A$ is corresponding to the pair of tableaux of shape $\lambda = ( \lambda_1, \lambda_2, \dots)$
  by the generalized RSK correspondence,
  then we have $\ell (w_A) = \lambda_1$.
\end{thm}
This theorem follows immediately from the following lemma.

\noindent
\begin{lem}
If $\pi= \alpha_1 \alpha_2 \dots \alpha_n \in \bP^n$ and 
$\alpha_k$ enters a marked tableau $S_{k-1}$ in the $j$th column  (of the first row),
then the longest increasing subsequence in $\pi$ ending in $\alpha_k$ has length $j$.
\end{lem}
\begin{proof}
We prove the claim by induction on $k$.
The result is trivial for $k=1$.
Suppose that it holds for $k-1$.

First we need to show the existence of an increasing subsequence of length $j$ ending in $\alpha_k$.
Let $\beta$ be the element of $S_{k-1}$ in the cell $(1,j-1)$.
Then we have $\beta < \alpha_k$, or $\beta$ and $\alpha_k$ are identical and unmarked,
since $\alpha_k$ enters in the $j$th column.
By induction, there is an increasing subsequence $\sigma$ of length $j-1$ ending in $\beta$.
Thus $\sigma \alpha_k$ is the desired subsequence.

Now we have to prove that there is no longer increasing subsequence ending in $\alpha_k$.
Suppose that such a sequence exists and let $\alpha_i$ be the preceding element of $\alpha_k$ in the subsequence.
Then it is satisfied that $\alpha_i < \alpha_k$, or $\alpha_i$ and $\alpha_k$ are identical and unmarked.
Since the sequence obtained by erasing $\alpha_k$ is a subsequence
whose length is greater than or equal to $j$ and whose ending is $\alpha_i$, 
by induction,
$\alpha_i$ enters in some $j'$th column such that $j' \geq j$ when $\alpha_i$ is inserted.
Thus the element $\gamma$ in the cell $(1,j)$ of $S_i$ satisfies
$\gamma \leq \alpha_i$,
so that 
$\gamma < \alpha_k$, or $\gamma$ and  $\alpha_k$ are identical and unmarked.

But since $\alpha_{k}$ is the element in the cell $(1,j)$ of $S_k$ and $i<k$,
we see that $\gamma > \alpha_k$, or that $\gamma$ and $\alpha_k$ are identical and marked.
It is a contradiction.
Therefore the lemma follows.
\end{proof}

For a given $\bP$-matrix $A$,
$\lambda_1$, where $\lambda= (\lambda_1, \lambda_2, \dots)$ is obtained from $A$ by the generalized RSK algorithm above,
gives the length of the longest increasing subsequence in a biword $w_A$.
On the other hands, 
$\lambda_1$ obtained by the shifted RSK algorithm between $\bP$-matrices and pairs of shifted marked tableaux
gives the length of
the longest \textit{ascent pair} for a biword $w_A$ associated with $A$
(see in \cite{TracyWidom200?}, \cite{HoffmanHumphreys}).

%
%
\section{Measures on $\bP$-matrices}
%
%

We give a combinatorial aspect of the $t$-Schur measure using the facts stated in the preceding section.
Let $x=(x_1,x_2, \dots)$ and $y=(y_1,y_2,\dots)$ be variables satisfying
$0 \leq x_i, y_j \leq 1$ for all $i,j$,
and $t \leq 0$.
Let $\bP_{m,n}$ denote the set of all $\bP$-matrices of size $m \times n$.
We abbreviate $\ell(w_A)$ to $\ell(A)$ for a $\bP$-matrix $A$.
Define a measure depending on a parameter $t$ as follows.
Assume that matrix elements $a_{ij}$ in $A$ are distributed independently with
the following distributions associated with parameters $x_iy_j$:
 \begin{align*}
  \Prob_t(a_{ij} = k ) &= \frac{1- x_i y_j}{1 - t x_i y_j} (x_i y_j)^k, \\    
  \Prob_t(a_{ij} = k') &= \frac{1- x_i y_j}{1 - t x_i y_j} (-t) (x_i y_j)^k
 \end{align*}
for $k \geq 1$ and
 $$
  \Prob_t(a_{ij} =0) = \frac{1- x_i y_j}{1 - t x_i y_j}.
 $$
This $\Prob_t$ indeed defines a probability measure on $\bP\cup \{ 0 \}$.
Actually we have
 $$ 
  \sum_{k=0}^{\infty} \Prob_t(|a_{ij}|= k) = 1
 $$
and $\Prob_t(a_{ij} = k') \geq 0$ for every $k$ since $t \leq 0$. 

Let 
 $$
   \bP_{m,n,s,u,r} := \{ A \in \bP_{m,n} | \ \sum_{1 \leq i \leq m} |a_{ij}| = s_j, 
    \sum_{1 \leq j \leq n} |a_{ij}| = u_i, 
    \mathrm{mark}(A)=r \}
 $$
for $s \in \bZ_{\geq 0}^n$, $u \in \bZ_{\geq 0}^m$, and $0 \leq r \leq mn$.
Here $\mathrm{mark}(A)$ is the number of marked entries in $A$. 
Then we have 
 \begin{equation} \label{uniform}
  \Prob_t(\{ A \}) = \prod_{\begin{subarray}{c} 1 \leq i \leq m \\ 1 \leq j \leq n \end{subarray} }
  \( \frac{1- x_i y_j}{1- t x_i y_j} \) (-t)^r x^s y^u 
  = \frac{1}{Z_t} (-t)^r x^s y^u
 \end{equation}
for $A \in \bP_{m,n,s,u,r}$.

If $\Prob_{t,m,n}$ denotes the probability measure obtained by putting $x_i=y_j=0$ for $i >m$ and $j>n$,
then it follows from \eqref{uniform}, Theorem 2, Theorem 3, and \eqref{generate} that
\begin{align*}
  \Prob_{t,m,n} (\ell \leq h) &= \Prob_t ( \{ A \in \bP_{m,n} \  | \ \ell(A) \leq h \} ) \\
                          &= \sum_{s,u,r} \Prob_t ( \{ A \in \bP_{m,n,s,u,r} \ | \ \ell(A) \leq h \} ) \\                       
                          &= \sum_{s,u,r} \# \{ A \in \bP_{m,n,s,u,r} \ | \ \ell (A) \leq h \} 
                            \  \frac{1}{Z_t} (-t)^r x^s y^u \\
                          &= \sum_{s,u,r} \# \{ (S,U) \ | \ \mathrm{type} \  s \  \mathrm{and} \  u, \ \mathrm{mark}(S)=r,\ 
                                   \lambda_1 \leq h \} \ \frac{1}{Z_t} (-t)^r x^s y^u \\
    &= \frac{1}{Z_t} \sum_{\Gessel} S_{\lambda}(x_1, \dots, x_m;t) s_{\lambda}(y_1, \dots, y_n).
\end{align*}
A set $\{ (S,U) | \ \mathrm{type} \  s \  \mathrm{and} \  u, \mathrm{mark}(S)=r, \ \lambda_1 \leq h \}$ 
consists of all pairs $(S,U)$ of the same shape $\lambda$ such that $\lambda_1 \leq h$,
where $S$ is a marked tableau which is of type $s$ and $\mathrm{mark}(S)=r$, 
and $U$ is a semi-standard tableau which is of type $u$.

Observe that the rightmost hand side in the above equality is the value of the $t$-Schur measure
with respect to a set $\{\lambda \in \P | \lambda_1 \leq h \}$.
Particularly, when $t=0$,
this measure on $\bP$-matrices turns to be the measure on $\bN$-matrices and 
it corresponds to the (original) Schur measure
(see \cite{Johansson1}, Johansson's $q$ is equal to our $\alpha^2$). 
On the other hand,
when $t=-1$, by the shifted RSK correspondence
we see that
it corresponds to 
the shifted Schur measure (see \cite{TracyWidom200?}).

\begin{remark}
The $t$-Schur measure at $t=-1$ does not coincide with the shifted Schur measure 
since the correspondences between $\bP$-matrices and partitions are different.
In fact,
Theorem 1 states that a (centered and normalized) limit distribution of  
$\ell(w_A)$ is identical with one of the length $L(w_A)$ of the longest ascent pair for $w_A$.
\end{remark}

%
%
\section{Proof of Theorem \ref{mainthm}}
%
%

In this section, 
we prove Theorem 1 using the methods developed in \cite{TracyWidom200?}.
We recall some notations.
Denote the Toeplitz matrix $T(\phi)=(\phi_{i-j})_{i,j \geq 0}$ and the Hankel matrix $H(\phi) = (\phi_{i+j+1})_{i,j \geq 0}$,
where $(\phi_n)_{n \in \bZ}$ is the sequence of Fourier coefficients of a function $\phi$
(see \cite{BottcherSilbermann}).
These matrices act on the Hilbert space $\ell^2 (\bZ_{+})$ $(\bZ_{+} = \bN \cup {0})$.
Also we put $T_h(\phi)= (\phi_{i-j})_{0 \leq i,j \leq h-1}$ and
$\Tilde{\phi} (z) := \phi(z^{-1})$. 
Let $P_h$ be the projection operator from $\ell^2(\bZ_{+})$ onto the subspace $\ell^2(\{0,1, \dots, h-1 \})$
and set $Q_h := I - P_h$, where $I$ is the identity operator on $\ell^2(\bZ_{+})$.

The following lemmas are keys in the proof.
The first one, Lemma 2 is a generalization of the Gessel identity 
(see \cite{Gessel}, \cite{TracyWidom2001}). 
\begin{lem}
\begin{equation} \label{Gessel}
    \sum_{\Gessel} S_{\lambda} (x;t) s_{\lambda} (y) = \det T_h (\Tilde{E}_{x,t} E_y)
 \end{equation}     
where $E_{x,t}$ is defined in \eqref{E} and $E_y = E_{y,0}$.
\end{lem}
\begin{proof}
Let $M(x;t)$ be the $\infty \times h$ submatrix $(e_{i-j}(x;t))_{i \geq 1, 1 \leq j \leq h}$ 
of the Toeplitz matrix $T (E_{x,t})$,
and for any subset $S \subset \bN$, 
let $M_S(x;t)$ be the submatrix of $M(x;t)$ obtained from rows indexed by elements of $S$.
In particular,  write $M(x)=M(x;0)$.

For a partition $\lambda = (\lambda_1, \lambda_2, \dots)$ such that $\lambda_1 \leq h$,
let $\lambda'$ (whose length is smaller than or equal to $h$) be the conjugate partition of $\lambda$
and let $S=\{ \lambda'_{h+1-i}+i | 1 \leq i \leq h \}$. 
Then we have
$\det M_S(x;t) = \det (e_{\lambda'_{h+1-i}+i-j}(x;t))_{1 \leq i,j \leq h}$.
Reversing the order of rows and columns in this determinant,
we have
$\det M_S(x;t) = \det (e_{\lambda'_i -i+j}(x;t)) = S_{\lambda}(x;t)$ by \eqref{S}.
In particular, $\det M_S (y) = \det ( e_{\lambda'_i -i +j}(y)) = s_{\lambda}(y)$.
It follows that 
 $$
  \sum_{\Gessel} S_{\lambda}(x;t) s_{\lambda}(y) =
  \sum_S \det M_S (x;t) \det M_S (y)
 $$
where the sum is all over $S \subset \bN$ such that $\# S= h$.
Then by the Cauchy-Binet identity, we have
  $$
   =  \det \  ^tM(x;t) M(y)  = \det T_h (\Tilde{E}_{x,t} E_y).
  $$
Hence the lemma follows.
\end{proof}

\begin{lem}
 If we put $\phi=\Tilde{E}_{x,t} E_y$, 
 then we have
   \begin{equation} \label{Fredholm}
      \det T_h(\phi) = E(\phi) \det( I-H_1 H_2) \mid_{ \ell^2 ( \{h,h+1,\dots\} ) }
   \end{equation} 
 where put $ H_1= H(\Tilde{E}_{x,t} E_y^{-1})$, $H_2= H (E_{x,t}^{-1} \Tilde{E}_y)$
 and $E(\phi):= \exp \{ \sum_{k=1}^{\infty} k (\log \phi)_k (\log \phi)_{-k} \}$.
 Here the determinant on the right side in \eqref{Fredholm} is a Fredholm determinant defined by
 $$    
  \det (I-K) \mid_{ \ell^2 ( \{h,h+1,\dots\} ) } := \det(Q_h - Q_h K Q_h)
 $$
 for any trace class operator $K$.
\end{lem}
\begin{proof}
We observe that both $H_1$ and $H_2$ are Hilbert-Schmidt operators.
We obtain the equality \eqref{Fredholm} by applying directly the relation 
between the Toeplitz determinant and the Fredholm determinant 
by Borodin and Okounkov \cite{BorodinOkounkov}, \cite{BasorWidom} to $\phi$.
We leave the detail to the reader.
\end{proof}

Note that for $\phi=\Tilde{E}_{x,t} E_y$, we have $E(\phi)=Z_t$.

If we denote by $J$ the diagonal matrix whose $i$-th entry equals $(-1)^i$,
then we have $\det (I-H_1 H_2) = \det (I- J H_1 H_2 J)$.
In general, we note that
$-J H (\phi(z)) J = H (\phi(-z))$ for any function $\phi(z)$. 

From \eqref{t-Schurmeasure}, \eqref{Gessel} and \eqref{Fredholm},
we obtain
 \begin{equation} \label{Schur-Fredholm}
    \PP(\lambda_1 \leq h) = \det (Q_h - Q_h J H_1 H_2 J Q_h ).
 \end{equation}

We make here $\alpha$-specialization and scaling.
We put $x_i=y_j= \alpha$ $(0<\alpha<1)$ 
for $1 \leq i \leq m$ and $1 \leq j \leq n$,
and $x_i=y_j=0$ for $i>m$ and $j>n$.
Put $\tau = m/n >0$.
We set $i= h+n^{1/3}x$ and $j=h+n^{1/3}y$, where $h= cn + n^{1/3} s$.
The positive constant $c$ will be determined later.

It is convenient to replace $\ell^{2}(\{h, h+1, \dots, \})$ by $\ell^2 (\bZ_{+})$.
Let $\Lambda$ be a shift operator on $\ell^2(\bZ_{+})$, 
i.e., $\Lambda \mathbf{e}_j = \mathbf{e}_{j-1}$ 
for the canonical basis $\{ \mathbf{e}_j \}_{j \geq 0}$ of $\ell ^2 (\bZ_{+})$
and $\Lambda^{*}$ the adjoint operator of $\Lambda$.
In $Q_h J H_1 H_2 J Q_h$, 
we interpret
a $Q_h$ appearing on the left 
as $\Lambda^h$ and 
a $Q_h$ appearing on the right as $\Lambda^{* h}$.
Then the $(i,j)$-entry of $ - Q_h J H_1 J $ is
 $$
    \frac{1}{2 \pi \sqrt{-1}} 
    \int \( \frac{z- \alpha}{z- t \alpha} \)^m \( \frac{1}{1- \alpha z} \)^n z^{-cn - n^{1/3} s -i -j} \frac{\d z}{z^2}
 $$
and the $(i,j)$-entry of $ - J H_2 J Q_h$ is
 $$ 
    \frac{1}{2 \pi \sqrt{-1}} 
    \int \( \frac{1- t \alpha z}{1- \alpha z} \)^m (1- \alpha z^{-1})^n z^{-cn - n^{1/3} s -i -j} \frac{\d z}{z^2},
 $$
where the contours of the integrals are both the unit circle.
For the second integral, 
if we change the variable $z \to z^{-1}$,
we obtain
 $$
    \frac{1}{2 \pi \sqrt{-1}} \int \( \frac{z- t \alpha}{z - \alpha} \)^m (1 - \alpha z)^n z^{cn + n^{1/3} s +i +j} \d z.
 $$
Putting 
 $$
    \psi(z) = \( \frac{z- \alpha}{z- t \alpha} \)^m \( \frac{1}{1- \alpha z} \)^n z^{-cn},
 $$
we find that two integrals are rewritten as
 \begin{align} 
    \frac{1}{2 \pi \sqrt{-1}} \int \psi(z) z^{- n^{1/3} s -i -j} \frac{\d z}{z^2}, \label{integral1} \\ 
    \frac{1}{2 \pi \sqrt{-1}} \int \psi(z)^{-1} z^{n^{1/3} s +i +j} \d z. \label{integral2}
 \end{align}

To estimate these integrals, we apply the steepest descent method.     
At first we determine the constant $c$.
Let $\sigma (z) = n^{-1} \log \psi (z)$, so that
 \begin{equation} \label{sigma'}
    \sigma' (z) =  \frac{\tau \alpha (1-t)}{(z- \alpha) (z -t \alpha)} + \frac{\alpha}{1- \alpha z} - \frac{c}{z}.
 \end{equation}
We choose a constant $c$ such that $\sigma(z)$ has the point $z$
satisfying the equality $\sigma' (z) = \sigma''(z) =0$. 
Then we obtain
 \begin{equation} \label{z_0}
    \frac{\tau (1-t)(z^2- t \alpha^2)}{(z - \alpha)^2 (z- t \alpha)^2} - \frac{1}{(1- \alpha z)^2} =0.
 \end{equation}
Since we assume $t \leq 0$,
it is immediate to see that
the function on the left hand side in \eqref{z_0} is strictly decreasing 
from $+\infty$ to $-\infty$ on the interval $(\alpha, \alpha^{-1})$ and
it follows that there is a unique point $z_0$ in $(\alpha, \alpha^{-1})$, where the left is equal to zero.
This is a saddle point and we set 
 \begin{equation} \label{c}
    c:= \alpha z_0 \( \frac{\tau (1-t)}{(z_0- \alpha) (z_0 - t \alpha)} + \frac{1}{1 - \alpha z_0} \)
 \end{equation}
from \eqref{sigma'}.
The constant $c$ is positive since $\alpha<z_0<\alpha^{-1}$ and $t \leq 0$.

It is clear that
the number counting with the multiplicity of zeros of the function $\sigma'(z)$ are three from \eqref{sigma'}
when $t \not=0$.
In the case where $t<0$, $\sigma' (z)$ has a zero in $(t \alpha, 0)$
because $\lim_{z \downarrow t \alpha} \sigma'(z) = - \infty$ and $\lim_{z \uparrow 0} \sigma'(z) = + \infty$. 
On the other hand, in the case where $t=0$,
the number of zeros of the function $\sigma'(z)$ are two.
Therefore, since $\sigma'(z)$ has a double zero $z_0$,
we have $\sigma'''(z_0) \not= 0$.
Further, since $\lim_{z \downarrow \alpha} \sigma' (z) = +\infty$ and $\lim_{z \uparrow \alpha^{-1}} \sigma' (z)= +\infty$, 
$\sigma'(z)$ is positive on $(\alpha, \alpha^{-1})$ except $z_0$,
so that $\sigma'''(z_0)$ is positive.

We call $\Gamma_{+}$ a steepest descent curve for the first integral \eqref{integral1} 
and $\Gamma_{-}$ for the second integral \eqref{integral2}.
On $\Gamma_{+}$,
the absolute value $|\psi (z) | = \exp \mathrm{Re}\  \sigma(z)$ is maximal at $z=z_0$ and
strictly decreasing as we move away from $z_0$ on the curve.
The first one emanates from $z_0$ at angles $\pm \pi /3$ with branches going to $\infty$ in two directions. 
The second one emanates from $z_0$ at angles $\pm  2\pi /3$ and is getting close to $z=0$.
 
Let $D$ be a diagonal matrix whose $i$-th entry is given by $\psi (z_0)^{-1} z_0^{n^{1/3}s +i}$, and multiply
$Q_h J H_1 H_2 J Q_h= (-Q_h J H_1 J) (-J H_2 J Q_h)$ by $D$ from the left and
by $D^{-1}$ from the right.
Note that the determinant of the right hand side in \eqref{Schur-Fredholm} is not affected.

By the discussion in Section 6.4.1 in \cite{TracyWidom200?}, 
we see that
the $([n^{1/3}x], [n^{1/3}y])$-entry of 
$$
  n^{\frac{1}{3}} D Q_h J H_1 H_2 J Q_h D^{-1}
$$
converges to $g K_{\mathrm{Airy}} ( g(s+x), g(s+y) )$ in the trace norm as $n \to \infty$, where we put
 \begin{equation} \label{g}
    g:=z_0^{-1} \(\frac{2}{\sigma'''(z_0)} \)^{1/3}
 \end{equation}
and $K_{\text{Airy}} (x,y)$ is the Airy kernel.
Clearly, $g$ is positive. 

Hence, by \eqref{Schur-Fredholm}, we have
\begin{align*}
   &  \lim_{n \to \infty} \PP_{\sigma} ( \lambda_1 \leq cn + n^{1/3} s) \\
 = &  \sum_{k=0}^{\infty} (-1)^k \sum_{0 \leq l_1< l_2< \dots < l_k} 
      \det \( ( D Q_h J H_1 H_2 J Q_h D^{-1})_{l_i, l_j} \)_{1 \leq i,j \leq k} \\
 = & \sum_{k=0}^{\infty} \frac{(-1)^k}{k !} \int_{[0, \infty)^k} 
      \det \( ( D Q_h J H_1 H_2 J Q_h D^{-1})_{[x_i], [x_j]} \)_{1 \leq i,j \leq k}
      \   \d x_1 \dots \d x_k \\
 = & \sum_{k=0}^{\infty} \frac{(-1)^k}{k !} \int_{ [0, \infty)^k} 
      \det \( n^{1/3} ( D Q_h J H_1 H_2 J Q_h D^{-1})_{[n^{1/3}x_i], [n^{1/3}x_j]} \)_{1 \leq i,j \leq k}
      \   \d x_1 \dots \d x_k \\
 \longrightarrow &  \sum_{k=0}^{\infty} \frac{(-1)^k}{k !} \int_{[0, \infty)^k} 
      \det \( g K_{\mathrm{Airy}}(g(s+x_i), g(s+x_j)) \)_{1 \leq i,j \leq k}
      \   \d x_1 \dots \d x_k \\
 = & \sum_{k=0}^{\infty} \frac{(-1)^k}{k !} \int_{[g s, \infty)^k} 
      \det \( K_{\mathrm{Airy}} (x_i, x_j)  \)_{1 \leq i,j \leq k}
      \   \d x_1 \dots \d x_k \\
 = & F_2 (gs).
\end{align*}

This shows the assertion of the theorem, 
where the constants $c_1 ( \alpha, \tau, t)$ and $c_2 (\alpha, \tau, t)$
are given by $c$ and $g^{-1}$, respectively.
We complete the proof of the theorem.

\begin{remark}
 In the case where $t=0$, we have $z_0 = \frac{ \alpha + \sqrt{\tau} }{1+ \sqrt{\tau} \alpha}$ by  \eqref{z_0}.
 Therefore we obtain
  $$
   c_1 (\alpha, \tau, 0) = \frac{(1+ \sqrt{\tau} \alpha)^2}{1- \alpha^2} -1 
  $$
 by \eqref{c} and
  $$
   c_2 (\alpha, \tau, 0) = g^{-1} 
   = \frac{ \alpha^{1/3} \tau^{-1/6} }{ 1-\alpha^2 } ( \alpha + \sqrt{\tau} )^{2/3} (1+ \sqrt{\tau} \alpha)^{2/3}
  $$
 by \eqref{g}. 
These values give the corresponding values in Theorem 1.2 in \cite{Johansson1}. 
Note that the relation between our $\alpha$ and Johansson's $q$ is given by $q = \alpha^2$.  
\end{remark}

\noindent
{\large\bf Acknowledgements}

I feel grateful to my supervisor, Professor Masato Wakayama, for suggestions and valuable comments. 
Also I would like to thank Professor Kazufumi Kimoto for detailed comments and encouragement.

%
%


\textsc{Sho Matsumoto}\\
Graduate School of Mathematics, Kyushu University.\\
Hakozaki Fukuoka 812-8581, Japan.\\
e-mail : \texttt{ma203029@math.kyushu-u.ac.jp}\\

\end{document}